\documentclass[10pt]{amsart}
\usepackage{amsmath,amssymb,enumerate}
\usepackage{hyperref}
\usepackage{enumitem}
\usepackage{breqn}
\newtheorem{thm}{Theorem}[section]
 \numberwithin{equation}{section} 
 \numberwithin{figure}{section} 
 \theoremstyle{plain}
 \theoremstyle{plain}    
 \newtheorem{cor}[thm]{Corollary} 
 \theoremstyle{plain}    
 \newtheorem{prop}[thm]{Proposition} 
 \newtheorem{defi}[thm]{Definition}
 \theoremstyle{plain}    
 \newtheorem{lem}[thm]{Lemma} 
 \theoremstyle{remark}
 
 \theoremstyle{definition}

\date{\today}

\author{Songchen Liu}
\address{School of mathematical sciences, Capital Normal University, Beijing, China}
\email{2230502073@cnu.edu.cn}

\begin{document}

\title{CHARACTERIZING THE RANGE OF THE COMPLEX MONGE--AMPÈRE OPERATOR}

\begin{abstract}
    In this note, we solve the complex Monge--Ampère equation for measures with a pluripolar part in compact Kähler manifolds. This result generalizes the classical results obtained by Cegrell in bounded hyperconvex domains. We also discuss the properties of the complex Monge--Ampère operator in some special cases.
\end{abstract}

\maketitle

\section{Introduction}
 Let $(X, \omega)$ be a Kähler manifold of complex dimension $n$, and let $\{\alpha\}\in H^{1,1}(X,\mathbb{R})$. The study of the complex Monge--Ampère type equation
\begin{equation}\label{eq 1.1}
(\alpha+dd^c \varphi)^n = \mu
\end{equation}
  has been one of the most important topics in Kähler geometry over the past few decades.\par
  
  When $\alpha>0$ and $\mu$ is a smooth volume form, Yau \cite{Yau78} proved that (\ref{eq 1.1}) admits a unique smooth solution, solving the famous Calabi conjecture. In this case, the left-hand side of (\ref{eq 1.1}) is the classic $n$-th wedge product.\par
  
 When the cohomology class $\{ \alpha\}$ is allowed to be degenerate and the measure 
$\mu$ is nonsmooth but non-pluripolar, the solution of (\ref{eq 1.1}) is closely related to the singular Kähler--Einstein metric.
To solve the corresponding equation, the concept of the non-pluripolar product plays a crucial role; see, for example, \cite{BEGZ10, DDL23, EGZ09, GZ17}.
In this setting, the left-hand side of (\ref{eq 1.1}) is given by the non-pluripolar product.\par
  
  When $\alpha = \omega$ and $\mu$ are allowed to have a pluripolar part (see Corollary \ref{cor 4.4}), many uncertainties arise.
  Coman, Guedj, and Zeriahi in \cite{CGZ08} developed theories to consider related problems. An important concept is the {\it domain of the definition of the complex Monge--Ampère operator} $DMA(X,\omega)$, see \cite[Page 394]{CGZ08}.\par
  
  A special and interesting problem mentioned in \cite[Question 12]{DGZ16} is that when $X= {\rm BL}_p \mathbb{P}^2,~p\in\mathbb{P}^2$, $0< \alpha:= \omega \in -c_1(K_X)$ and $\mu : = [D]\wedge \omega/V$, where $D$ is a smooth anticanonical divisor, $[D]$ is the (1,1)-current integral along $D$ and $V$ is the normalized constant, can one find $\varphi \in DMA(X,\omega)$ such that ${\rm MA}_\omega(\varphi) = \mu$?\par
  
  Motivated by these related problems, we study the complex Monge–Ampère equation with a measure that has a pluripolar part in this note. Our focus will be on the 
  Błocki--Cegrell class $\mathcal{D}(X,\omega) \subset DMA(X,\omega)$ (see Definition \ref{def 2.10}).\par
  
   In the recent paper \cite{AWNW22}, Andersson, Witt Nyström, and Wulcan introduced the finite non-pluripolar energy class $\mathcal{G}(X,\omega)$. 
  They proved that $\mathcal{D}(X,\omega)\subset \mathcal{G}(X,\omega)$; i.e., if $\varphi \in \mathcal{D}(X,\omega)$, then $\varphi \in L^1 \Big(\langle \omega_\varphi^i \rangle\wedge \omega^{n-i} \Big),~i=0,1,...,n-1$, see \cite[Theorem 1.10]{AWNW22}.\par
  
  In Section \ref{sec 3}, we show that one benefit of \cite[Theorem 1.10]{AWNW22} is that the Monge--Ampère operator on $\mathcal{D}(X,\omega)$ can be expressed via the non-pluripolar product. Furthermore, we present the corresponding generalized theorem.
 \begin{thm}$(=${\rm Theorem} \ref{thm 3.2}$)$ 
     Let $(X,\omega)$ be a compact Kähler manifold of complex dimension $n$, and let $\varphi_i \in \mathcal{D}(X,\omega),~ i=1,...,n$. Then, $\varphi_1 \in L^1\Big(\langle (\omega+dd^c \varphi_2) \wedge...\wedge (\omega+ dd^c \varphi_n) \rangle\wedge\omega \Big)$.\par
     If $X$ is a compact Kähler surface and $\varphi,\psi \in \mathcal{D}(X,\omega)$, then we have
    $$
    {\rm MA}_\omega(\varphi,\psi)= \frac{1}{2}\left[ (\omega+dd^c\varphi)\wedge\langle \omega+dd^c \psi \rangle +(\omega+dd^c\psi)\wedge\langle \omega+dd^c \varphi \rangle \right].
    $$
    
     {\rm The~definition of 
     ${\rm MA}_\omega(\cdot,...,\cdot)$ can be found in Proposition~\ref{prop 2.13}}.
 \end{thm}

  %Now, we consider the complex Monge--Ampère equation with a normalized measure $\mu$ on $X$.
 Now, we consider the equation $ {\rm MA}_\omega(\varphi) = \mu$, where  $\varphi \in \mathcal{D}(X,\omega)$ and $\mu$ is a normalized positive Borel measure.
  The Cegrell--Lebesgue decomposition, Corollary \ref{cor 4.4}, indicates that we have the unique decomposition $\mu = \mu_r + \mu_s$, where $\mu_r$ is the non-pluripolar measure and $\mu_s = \mathbf{1}_{\{u=-\infty\}}\mu$ for some $u\in {\rm PSH}(X,\omega)$. 
   A natural approach is to decompose ${\rm MA}_\omega(\varphi)$, then cosinder ${\rm MA}_\omega (\varphi)_r = \mu_r$, ${\rm MA}_\omega(\varphi)_s = \mu_s$.\par

 Thanks to \cite[Corollary 1.8]{CGZ08}, which yields ${\rm MA}_\omega(\varphi)_r = \langle \omega_\varphi^n \rangle$ and ${\rm MA}_\omega(\varphi)_s = \mathbf{1}_{\{ \varphi = -\infty \}} {\rm MA}_\omega(\varphi)$. 
 Consequently, ${\rm MA}_\omega(\varphi) = \mu$ reduces to the following equations:
 $$
 \langle \omega^n_\varphi \rangle  = \mu_r~~{\rm and}~~\mathbf{1}_{\{\varphi = -\infty\}}{\rm MA}_\omega(\varphi) = \mu_s,~\varphi \in \mathcal{D}(X,\omega).
 $$
 Here the complex Monge--Ampère equation on the left-hand side, defined using the non-pluripolar product, is the familiar one.
 Thus, by combining it with the recent work of Darvas, Di Nezza, and Lu, \cite{DDL23}, named {\it relative pluripotential theory}, we can solve the complex Monge--Ampère equation in a particular case.
\begin{thm}$(=${\rm Theorem} \ref{thm 4.9}$)$
    Let $(X,\omega)$ be a compact Kähler manifold of complex dimension $n$ such that $\int_X \omega^n =1$, and let $\mu_s$ be a pluripolar measure on $X$ supported on some pluripolar set. If there exists $\varphi \in \mathcal{D}(X,\omega)$ such that 
    ${\rm MA}_\omega(\varphi)_s = \mu_s$, then for all $0 \leq f \in L^p(\omega^n),~p>1$ such that $\int_X f \omega^n = 1 - \int_X \mu_s$, there exists a $\psi \in {\rm PSH}(X,\omega)$ that is the solution of the equation
    $$
    {\rm MA}_\omega(\psi)= f \omega^n + \mu_s,~\psi \in \mathcal{D}(X,\omega).
    $$
\end{thm}

\subsection{Comparing the classical result in the local setting} Let $\Omega \subset   \mathbb{C}^n$ be a bounded hyperconvex domain. In \cite{ACCP09,BT76,Ce98,Ce04}, a powerful theory for the Dirichlet problem of the complex Monge--Ampère equation on $\Omega$ was developed by several authors. In particular, the author was inspired by \cite{Ce98,Ce04} to consider the corresponding global result. Let $\mu$ be a positive Radon measure on $\Omega$. Cegrell proved the following result.
\begin{thm}$($\cite[Theorem 6.2]{Ce04}$)$
 Assume that $\mu = f(dd^c \varphi_0)^n + \mu_s$, where $\varphi_0 \in \mathcal{E}_0(\Omega)$ and $f\in L^1\left( (dd^c \varphi_0)^n \right)$. If there exists $\psi \in \mathcal{F}(\Omega)$ such that ${\rm MA}(\psi)= \mu_s$, then there exists $\varphi \in \mathcal{F}(\Omega)$ such that
$$
{\rm MA}(\varphi) = \mu.
$$

{\rm The definition of $\mathcal{F}$ and $\mathcal{E}_0$ can be found in \cite{Ce04}.}
\end{thm}
 In the global setting, we replace the condition of the pluripolar part of $\mu$ with the condition that there exists $\varphi \in \mathcal{D}(X,\omega)$ such that 
 ${\rm MA}_\omega(\varphi)_s = \mu_s$ and require that $\mu_r = f\omega^n$ for some $0\leq f\in L^p(\omega^n),~p>1$.

\subsection{Application} Let us focus on compact Kähler surfaces $(X,\omega)$. Assume that $\int_X\omega^2=1$, and let $\mu$ be a positive Borel measure on $X$, supported on a pluripolar set. Suppose there exists $\varphi \in \mathcal{D}(X,\omega)$ such that ${\rm MA}_\omega(\varphi) = \mu$. By applying the two main theorems of this note, we obtain the following result.
\begin{thm}
    For any $0\leq t \leq 1$ and $0\leq f \in L^p(\omega^2),p>1$ such that $\int_X f \omega^2 =1$, there exists $\psi_{f,t} \in \mathcal{D}(X,\omega)$ that satisfies
    $$
    {\rm MA}_\omega(\psi_{f,t}) = (1-t) f\omega^2 + t\mu.
    $$
\end{thm}

\subsection{Notation and Conventions}
In this note, unless stated otherwise, we always assume that
\begin{itemize}[label=\textbullet, font=\Large]
\item  In the local setting, we consider the domain $\Omega$ in $\mathbb{C}^n$.  In the global setting, we consider the compact Kähler manifolds $(X,\omega)$.
\item The operator $dd^c : = \frac{\sqrt{-1}}{\pi} \partial\bar{\partial}$ is the normalized $\sqrt{-1}\partial\bar{\partial}$ operator.
\item In the global setting, let ${\rm PSH}(X, \omega)$ be the set of upper semicontinuous functions $\varphi$ such that $dd^c \varphi + \omega \geq 0$. Plurisubharmonic functions are abbreviated as psh functions. When $\varphi \in {\rm PSH}(X,\omega)$, we abbreviate this by saying that $\varphi$ is a $\omega$-psh function.
\item In any setting, we always assume that a measure is a positive Borel measure.
\end{itemize}

\section{Preliminaries}
  In pluripotential theory, the definition of the complex Monge-Ampère operator is a central issue. In the local setting, let $u_1,...,u_k$ be a locally bounded psh function. Following the construction of Bedford--Taylor \cite{BT76} (see also \cite[Chapter 3]{GZ17}), where the current $dd^c u_1 \wedge...\wedge dd^c u_k$ is always well defined.

\subsection{Non-pluripolar products}
    Let $(X,\omega)$ be a compact Kähler manifold. In the global setting, given the positive (1,1)-currents $\omega^1 + dd^c\varphi_1, ..., \omega^p + dd^c\varphi_p$ on $X$. Assume $\omega^i$, $i=1,...,p$ are Kähler forms and $\varphi_i \in {\rm PSH}(X,\omega^i)$. Following the construction of Bedford--Taylor \cite{BT76} in the local setting, it was shown in \cite{BEGZ10} that the $d$-closed currents, so-called {\it non-pluripolar product}
$$
\langle \omega^1_{\varphi_1}\wedge...\wedge \omega^p_{\varphi_p} \rangle,
$$
is well-defined. By \cite[Proposition 1.4]{BEGZ10}, the non-pluripolar product is symmetric and multilinear.
For a $\omega$-psh function $\varphi$, the {\it non-pluripolar complex Monge--Ampère measure} of $\varphi$ is  $\langle \omega_\varphi^n \rangle:= \langle (\omega_\varphi)^n \rangle$.

 Let $Y$ be a complex manifold. A set $A \subset Y$ is said to be {\it pluripolar} in $X$ if, for all $z \in A $, there exists a holomorphic coordinate chart $U$ near $z$ in $Y$ and a psh function $\varphi \in {\rm PSH}(U)$ such that $A \cap U \subset \{  \varphi = -\infty\}$. A measure $\mu$ of $Y$ is said to be a {\it non-pluripolar measure} if $\int_A \mu=0$ for all pluripolar sets $A\subset Y$.
 In particular, when $p=n$, the measure $\langle \omega^1_{\varphi_1}\wedge...\wedge \omega^n_{\varphi_n} \rangle$ mentioned above is non-pluripolar.\par

  If $\varphi$ and $\psi$ are two $\omega$-psh functions on $X$, then $\psi$ is said to be {\it less singular} than $\varphi$, i.e., $\varphi \preceq  \psi$, 
 if they satisfy $\varphi \leq \psi + C$ for some $C \in \mathbb{R}$. We say that $\varphi$ has the {\it same singularity} as $\psi$, i.e., $\varphi\backsimeq  \psi$,
 if $\varphi\preceq \psi$ and $\psi \preceq\varphi$.
 Then, we have the Volume comparison, due to Witt Nyström \cite{WN19}.
\begin{thm}\label{thm 2.2}
Let $u,v \in {\rm PSH}(X,\omega)$ such that $v \preceq u$. Then,
$$
\int_X \langle \omega_{u}^n \rangle \geq
\int_X \langle \omega_{v}^n \rangle
$$
\end{thm}
%Note that in \cite{DDL18}, Darvas--Di Nezza--Lu proves a more general version of Theorem \ref{thm 2.2}. \par
%Now, we introduce two important concepts in pluripotential theory. \subsubsection{Envelopes} Let $f$ be a function on $X$ such that $f: X \to \mathbb{R}\cup \{ \infty \}$. We define the envelope of $f$ in the class ${\rm PSH}(X, \omega)$ as $$ P_\omega(f):= \big( \sup\{ \varphi \in {\rm PSH}(X,\omega): \varphi \leq f \} \big)^*, $$ where $h^*(x):= \limsup_{y \to x} h(y)$. In the particular case where $f : = \min(\psi, \phi)$, we denote the envelope as  $P_\omega(\psi, \phi) := P_\omega(\min(\psi, \phi))$.\par
 In our study of the complex Monge--Ampère equation, the following envelope construction is essential:
\begin{defi}
    Given $\psi,\varphi \in {\rm PSH}(X, \omega)$, the envelope with respect to singularity type $P_\omega[\psi]$ is defined by
  $$
  P_\omega[\psi]:= \left( \lim_{C\to +\infty} \sup \{ u \in {\rm PSH}(X,\omega):~u\leq \psi +C, u \leq 0 \}^* \right)^*.
  $$
\end{defi}
We summarize \cite[Remark 3.4,Theorem 3.14]{DDL23} as follows.
\begin{prop}\label{DDL23 2.4}
    Let $\varphi \in {\rm PSH}(X,\omega)$. Then, we have\\
    $(i)$. $\int_X \langle \omega_\varphi^n \rangle = \int_X \langle \omega_{P_\omega[\varphi]}^n \rangle$,\\
    $(ii)$. Set 
    $
    F_\varphi := \{ v \in {\rm PSH}^-(X,\omega):~\varphi \preceq v \leq 0 ~{\rm and}~\int_X \langle \omega_v^n \rangle = \int_X \langle \omega_\varphi^n \rangle \}
    $. One has $P_\omega[\varphi] = \sup F_\varphi$. In particular, $P_\omega[\varphi] \succeq \varphi$ and $P_\omega[\varphi]=P_\omega[P_\omega[\varphi]]$.
\end{prop}
A $\omega$-psh function $\varphi$ is said to be a {\it model potential} if $P_\omega[\varphi] = \varphi$ and $\int_X \langle \omega_\varphi^n \rangle>0$. 
Clearly, $P_\omega[\varphi]$ is a model potential for all $\varphi \in {\rm PSH}(X,\omega)$ so that $\int_X \langle \omega_\varphi^n \rangle>0$.
%\subsubsection{Relative full mass class} 
\begin{defi} 
Given a potential $\phi\in {\rm PSH}(X, \omega)$ such that $\int_X \omega_\phi^n>0$, the relative full mass class is defined by 
$$ \mathcal{E}(X, \omega, \phi):=\{ u \in {\rm PSH}(X,\omega): u \preceq \phi,~\int_X \langle  \omega_u^n \rangle= \int_X \langle \omega_\phi^n \rangle \}. 
$$
\end{defi}
For the case in which $\phi$ is a model potential, Darvas, Di Nezza, and Lu \cite{DDL23} achieved a series of significant results on pluripotential theory and the (non-pluripolar product) complex Monge--Ampère equation in the relative full mass class $\mathcal{E}(X,\omega,\phi)$\footnote{In their results, it is not necessary for $\omega$ to be Kähler. They only require $\{ \omega \}$ to be big.}. 

%Their results play a key role in our article. \begin{thm}\label{DDL23 2.6} $($\cite[Theorem 5.17]{DDL23}$)$ Let $\phi \in {\rm PSH}(X,\omega)$ be a model potential and $\int_X \langle \omega_\phi^n \rangle >0$. Assume that $\mu$ is a non-pluripolar positive measure such that $\mu(X) = \int_X \langle \omega_\phi^n \rangle$. Then, there exists a unique normalized $u \in \mathcal{E}(X,\omega,\phi)$ such that $\langle \omega_u^n\rangle = \mu $. \end{thm} \begin{thm}\label{DDL23 2.7}  $($\cite[Theorem 5.20]{DDL23}$)$   Let $\phi \in {\rm PSH}(X,\omega)$ be a model potential and $\int_X \langle \omega_\phi^n \rangle >0$. Let $u \in \mathcal{E} ( X , \omega , \phi )$ with $\sup_X u = 0$. If $\langle \omega_u^n \rangle = \mu$, where $\mu$ is a positive measure such that $\mu = f\omega^n,~f\in L^p(\omega^n)$ for some $p>1$, then $u$ has the same singularity type as $\phi$.  \end{thm}

\subsection{The complex Monge--Ampère operators on compact Kähler manifolds}\label{sec 2.2} %In this section, we assume that $(X,\omega)$ is a compact Kähler manifold such that $\int_X \omega^n >0$. \begin{defi}\label{def 2.8} Fix $\varphi \in {\rm PSH}(X,\omega)$. We say that the complex Monge--Ampère measure $(\omega+dd^c\varphi)^n$ is well defined and write $\varphi \in DMA(X,\omega)$ if there exists a Radon measure $\mu$ such that for any sequence $\{ \varphi_j \}$ of bounded $\omega$-psh functions decreasing to $\varphi$ on X, the Monge--Ampère measures $(\omega+dd^c \varphi_j)^n$ converge weakly to $\mu$. We set $$ {\rm MA}_\omega(\varphi)=(\omega + dd^c\varphi)^n := \mu. $$ {\rm It is easy to see that ${\rm PSH}(X,\omega)\cap L^\infty(X)\subset DMA(X,\omega)$, \cite[Chapter 3]{GZ17}. Moreover, when $\varphi \in {\rm PSH}(X,\omega)\cap L^\infty(X)$, ${\rm MA}_\omega(\varphi)$ is actually the Bedford--Taylor wedge product.} \end{defi}

  In the local setting, assume that $\Omega\subset \mathbb{C}^n$ open. 
  We denote $\mathcal{D}(\Omega)\subset {\rm PSH}^-(\Omega)$ as the largest subclass of the class of negative psh functions in which the complex Monge--Ampère operator can be well-defined, named the Blocki--Cegrell class; i.e., if $\varphi \in \mathcal{D}(\Omega)$, then there exists a measure $\mu$ such that for any sequence 
  $\{\varphi_j\} \subset {\rm PSH}\cap L^\infty_{loc}(\Omega)$, $\varphi_j \searrow \varphi$, the sequence $(dd^c \varphi_j )^n $ converges weakly to $\mu$.
  We denote ${\rm MA}(\varphi)= ( dd^c \varphi)^n := \mu$. In \cite{Blo04,Blo06}, Blocki carefully studied $\mathcal{D}(\Omega)$. When $\Omega$ is a bounded hyperconvex domain, Blocki \cite[Theorem 2.4]{Blo06} showed that $\mathcal{D}(\Omega)$ coincides with the class $\mathcal{E}(\Omega) \subset {\rm PSH}^-(\Omega)$ defined by Cegrell in \cite[Section 4]{Ce04}.\par

   %The bounded hyperconvex domain $\Omega \subset \mathbb{C}^n$ is a bounded domain such that there exists a bounded psh function $h :\Omega \to (-\infty,0)$ that satisfies $\{ h < c \}\Subset \Omega$ for all $c<0$. The unit ball in $\mathbb{C}^n$  is clearly a bounded hyperconvex domain.\par When $\Omega$ is a bounded hyperconvex domain, Blocki \cite[Theorem 2.4]{Blo06} showed that $\mathcal{D}(\Omega)$ coincides with the class $\mathcal{E}(\Omega) \subset {\rm PSH}^-(\Omega)$ defined by Cegrell in \cite{Ce04}. That is, if $0 \geq \varphi \in \mathcal{E}(\Omega)$, then, for $\forall U\Subset \Omega$, there exists a decreasing sequence $\varphi_j \in \mathcal{E}_0(\Omega)$ such that $\varphi_j \searrow \varphi$ in $U$ and $ \sup_j \int_\Omega (dd^c \varphi_j)^n <+\infty$. Where $\mathcal{E}_0(\Omega):= \{ \varphi\in {\rm PSH}^-\cap L^\infty(\Omega): \lim_{z\to \zeta} \varphi (z) =0 ,\forall \zeta \in \partial\Omega ~{\rm and} ~\int_\Omega (dd^c\varphi)^n <+\infty\}$.\par
 Now, we can consider the Blocki--Cegrell class in the global setting.
\begin{defi}\label{def 2.10} 
Let $\mathcal{D}(X,\omega)$ be the set of functions $\varphi \in {\rm PSH}(X,\omega)$ such that locally, on any small open holomorphic coordinate chart $U \subset X$, the psh function $\varphi|_U + \rho_U \in \mathcal{D}(U)$, where $\rho_U$ is the local potential of $\omega$ on $U$ such that $\varphi|_U + \rho_U \leq 0$.
\end{defi}
  It is easy to see that we have $\mathcal{D}(X,\omega) \subset DMA(X,\omega)$. Coman, Guedj, and Zeriahi characterized $\mathcal{D}(X,\omega)$ in \cite[Section 3]{CGZ08}. They found that $\mathcal{D}(X,\omega)$ is much smaller than $DMA(X,\omega)$.\par
  %Now we have \begin{table}[h] \centering \caption{Relations between several classes}\label{tab 1} \begin{tabular}{|c|c|} \hline $\Omega$ is a bounded domain & $L^\infty_{loc}\cap {\rm PSH}^-(\Omega) \subset \mathcal{D}(\Omega)$ \\\hline $\Omega$ is a bounded hyperconvex domain & $ \mathcal{E}_0(\Omega)\subset L^\infty_{loc}\cap {\rm PSH}^-(\Omega) \subset  \mathcal{D}(\Omega) = \mathcal{E}(\Omega)$ \\\hline $X$ is a compact Kähler manifold & $L^\infty(X) \cap {\rm PSH}(X,\omega) \subset \mathcal{D}(X,\omega)\subset DMA(X,\omega)$ \\ \hline \end{tabular}  \end{table}\par
  
  Let $\varphi \in \mathcal{D}(X,\omega)$ and $U\subset X$ be a holomorphic coordinate such that $\omega = dd^c g$ on $U$.
  Let $\chi \in C(X)$ such that ${\rm supp}\chi \Subset U$. 
  Set $\varphi_k := \max\{ \varphi, -k \}$; then, we have
\begin{align*}
    \int_X \chi {\rm MA}_\omega(\varphi) &= \lim_k \int_X \chi (\omega+dd^c \varphi_k)^n\\
    &= \lim_k \int_U \chi (dd^c (g + \varphi_k) )^n =  \int_U \chi {\rm MA}(\varphi+g)
\end{align*}
by the definition of the complex Monge--Ampère operator. Hence,
\begin{equation}\label{eq 2.10}
    \chi {\rm MA}_\omega(\varphi) = \chi \mathbf{1}_U {\rm MA}(\varphi+ g),
\end{equation}
where $\mathbf{1}_U MA(\cdot)$ represents the complex Monge--Ampère operator in $U$. We will use this concept frequently and will not explain it in detail again.\par

By \cite[Theorem 1.2]{Blo06}, one has
\begin{lem}\label{lem 2.11}
Let $\varphi \in \mathcal{D}(X,\omega)$ and $\psi \in {\rm PSH}(X,\omega)$ such that $\psi \succeq \varphi$. Then, $\psi \in \mathcal{D}(X,\omega)$.
\end{lem}
And the following is derived from \cite{Blo06,Ce98}:
\begin{lem}
    Let $\varphi, \psi \in \mathcal{D}(X,\omega)$. Then, $(1-t)\varphi+t \psi \in \mathcal{D}(X,\omega)$ for $~0\leq t \leq 1$, which means that $\mathcal{D}(X,\omega)$ is a convex set.
\end{lem}

Set ${\rm MA}_\omega(\varphi_1,...,\varphi_n):=(\omega+dd^c \varphi_1)\wedge...\wedge(\omega+dd^c \varphi_n),~\varphi_i \in {\rm PSH}(X,\omega)\cap L^\infty(X)$. We then have the following property:
\begin{prop}\label{prop 2.13}
    Let $\varphi_1,...,\varphi_n \in \mathcal{D}(X,\omega)$ and $\varphi^j_i \in {\rm PSH}(X,\omega)\cap L^\infty(X)$ such that $\varphi^j_i \searrow \varphi_i$ as $j\to +\infty$; then, ${\rm MA}_\omega(\varphi^j_1,...,\varphi^j_n)$ converges weakly to a Radon measure $\mu$, and the limit measure does not depend on the particular sequence. We denote $\mu$ by ${\rm MA}_\omega(\varphi_1,...,\varphi_n)$.
\end{prop}
\begin{proof}
     The similar conclusion in the bounded hyperconvex domain  $\Omega$ is a direct consequence of \cite[Theorem 4.2]{Ce04}. 
     The remaining work can be completed by finding a suitable coordinate and utilizing (\ref{eq 2.10}).
\end{proof}
  Let $\varphi_1,...,\varphi_n \in \mathcal{D}(X,\omega)$. We define
$ {\rm MA}_\omega\left({\varphi_1}^{(k_1)},...,{\varphi_n}^{(k_n)}\right)$
as ${\rm MA}_\omega(\psi_1,...,\psi_n)$, where $\psi_1=...=\psi_{k_1}:= \varphi_1$, $\psi_{k_1+1}=...=\psi_{k_1+k_2} := \varphi_2$,...,$\psi_n=\psi_{n-1}=...=\psi_{n-k_n}=\varphi_n$. 
\begin{cor}\label{cor 2.14}
    Assume that $\varphi_1,...,\varphi_s \in \mathcal{D}(X,\omega)$. Then, for $0\leq t_1,...,t_s \leq 1$ such that $\sum_i t_i =1$, we have $ \sum_i t_i \varphi_i  \in \mathcal{D}(X,\omega)$ and 
    \begin{align*}
    {\rm MA}_\omega &\left(\sum_i t_i \varphi_i\right)\\
    &= \sum_{k_i \geq 0, k_1+...+k_s=n} C_{k_1,...,k_s}\cdot t_1^{k_1}\cdot...\cdot t_s^{k_s} \cdot {\rm MA}_\omega\left(\varphi_1^{(k_1)},...,\varphi_s^{(k_s)} \right),
    \end{align*}
    where $C_{k_1,...,k_s} : = \frac{n!}{k_1!\cdot...\cdot k_s!}$.
\end{cor}
\begin{proof}
   By definition, it suffices to prove this result for bounded $\omega$-psh functions, which follows directly from the multilinear of Monge--Ampère operator.
\end{proof}
%\begin{rem}\label{rem 1} Let $\Omega \subset \mathbb{C}^n$ be a bounded hyperconvex domain and $u_1,...,u_n \in \mathcal{D}(\Omega)$. In the proof of Proposition \ref{prop 2.13}, the measure ${\rm MA}(u_1, ..., u_n)$ is well defined as a limiting measure. Then, we can define ${\rm MA}\left(u_1^{(k_1)},...,u_s^{(k_s)} \right)$ in the local setting, which is similar to the way it is defined in the global setting described in (\ref{eq 2.2}). We can then obtain the same conclusion in the local setting. That is, for $u_1,...,u_s \in \mathcal{D}(\Omega)$ and $0\leq t_1,...,t_s$, we have $ \sum_i t_i u_i  \in \mathcal{D}(\Omega)$ and   \begin{align*} {\rm MA} &\left(\sum_i t_i u_i\right)\\  &= \sum_{k_i \geq 0, k_1+...+k_s=n} C_{k_1,...,k_s}\cdot t_1^{k_1}\cdot...\cdot t_s^{k_s} \cdot {\rm MA}\left(u_1^{(k_1)},...,u_s^{(k_s)} \right)   \end{align*} where $C_{k_1,...,k_s} : = \frac{n!}{k_1!\cdot...\cdot k_s!}$. \end{rem}

\subsection{Non-pluripolar energy} 
Let us now introduce the non-pluripolar energies.
\begin{defi}\label{def 2.15}
We define the non-pluripolar energy of $\varphi \in {\rm PSH}(X,\omega)$ as
$$
E^{np}(\varphi):=\frac{1}{k+1}\sum_{j=0}^{n-1} \int_X \varphi \langle (dd^c \varphi + \omega)^j \rangle\wedge \omega^{n-j}.
$$
Then, we can define the finite energy class $\mathcal{G}(X,\omega):= \{ \varphi \in {\rm PSH}(X,\omega): E^{np}(\varphi) > -\infty \}$.
\end{defi}
Note that $\varphi \in \mathcal{G}(X,\omega)$ {\it iff} $\varphi \in L^1 \left(\langle  \omega_\varphi^k \rangle\wedge \omega^{n-k} \right),~k=0,1,...,n-1$ since $\varphi$ is bounded from above. Hence, if $\varphi \in \mathcal{G}(X,\omega)$, we can define the currents 
$$
[\omega+dd^c \varphi]^p:= (\omega+dd^c \varphi)\wedge \langle \omega_\varphi^{p-1} \rangle ~~{\rm and}~~S^{\omega}_p(\varphi):= [\omega+dd^c \varphi]^p - \langle \omega_\varphi^p \rangle,
$$
for $p=1,2,...,n$. In \cite[Proposition 3.3]{AWNW22}, they showed that $[\omega+dd^c \varphi]^p$ and $S^{\omega}_p(\varphi)$ are well-defined globally closed positive currents on $X$ that only depend on the current $dd^c \varphi+\omega$ and not on the choice of $\omega$ as a Kähler representative in the class $[\omega]$.\par

There are two useful results follow from \cite{AWNW22}:
\begin{thm}\label{Theorem 2.16}
$($\cite[Theorem 6.8]{AWNW22}$)$ 
    Assume that $\varphi \in \mathcal{G}(X,\omega)$ and that $\eta$ is a Kähler form in $[\omega]$ so that $\eta=\omega+dd^c g$, where $g $ is a smooth function on $X$. Let 
    $\varphi_l:=\max \{  \varphi,g-l \}$. Then, for $1 \leq p \leq n $, one has
$$
  \left( \omega+dd^c \varphi_l \right)^p \to [\omega+dd^c \varphi]^p + \sum_{j=1}^{p-1} S^\omega_j(\varphi)\wedge \eta^{p-j},~l\to \infty.
$$
\end{thm}

\begin{thm}\label{thm 2.17}$($\cite[Theorem 1.11]{AWNW22}$)$
 Let $(X,\omega)$ be a compact Kähler manifold of dimension $n$. Then,
 $$
 \mathcal{D}(X,\omega)\subset \mathcal{G}(X,\omega).
 $$
\end{thm}
In particular, $\mathcal{D}(X,\omega)\subset DMA(X,\omega)\cap \mathcal{G}(X,\omega).$

\section{The class of $DMA(X,\omega)\cap \mathcal{G}(X,\omega)$}\label{sec 3}
In this section, we assume that $(X,\omega)$ is a compact Kähler manifold of complex dimension $n$. Note that if $\varphi \in DMA(X,\omega)\cap \mathcal{G}(X,\omega)$, by the definition of the complex Monge--Ampère operator and Theorem \ref{Theorem 2.16}, we have
\begin{equation}\label{eq 3.1}
\begin{aligned}
{\rm MA}_\omega(\varphi) =  [\omega+dd^c \varphi]^n + \sum_{j=1}^{n-1} S^\omega_j (\varphi) \wedge \omega^{n-j}= \langle  \omega_\varphi^n \rangle +\sum_{j=1}^n S^\omega_j (\varphi) \wedge \omega^{n-j}.
\end{aligned}
\end{equation}
  It may be helpful to consider the solution of the complex Monge--Ampère equation 
  ${\rm MA}_\omega(\varphi)= \mu$ when $\mu$ is allowed to be a normalized measure with a pluripolar part. A similar idea was proposed in \cite[Remark 11.13]{AWNW22}, but they did not consider terms with $S^\omega_j(\varphi)\wedge\omega^{n-j},~j=1,...,n-1$.\par
  In particular, when $(X,\omega)$ is a compact Kähler surface, the above argument becomes more interesting. Let $f\in C^\infty(X)$ and $C\in \mathbb{R}^+$ so that $\eta:= C^{-1} dd^c f + \omega>0$. Set $g:= C^{-1} {f}$ and $\varphi_l':= \max\{ \varphi, g-l \}$, $\varphi_l :=\max \{ \varphi, -l \}$. Since $\varphi_l',\varphi_l\in {\rm PSH}(X,\omega)\cap L^\infty(X)$ and $\varphi_l',~\varphi_l$ decrease to $\varphi$ as $l\to \infty$, it follows from Theorem \ref{Theorem 2.16} and the definition of the complex Monge--Ampère operator that
\begin{equation}
\begin{split}
     {\rm MA}_\omega (\varphi)&=[\omega+dd^c \varphi]^2 + S^\omega_1(\varphi)\wedge\omega\\
                &= [\omega+dd^c \varphi]^2 + S^\omega_1(\varphi)\wedge(\omega+dd^c g).
\end{split}
\label{eq 3.2}
\end{equation}
Hence, for $ w\in C^\infty(X)$, it follows from (\ref{eq 3.2}) that
$$
\int_X w dd^c g \wedge S^\omega_1(\varphi)=\frac{1}{C} \int_X w dd^c f \wedge S^\omega_1=0.
$$
When $w=f$, we obtain $\int_X \left(\sqrt{-1}\partial f \wedge \bar{\partial} f \right) \wedge S^\omega_1(\varphi)=0,~\forall f \in C^\infty(X)$. Then, we can construct a positive (1,1)-form $\alpha$ such that $\alpha>\omega$ and 
$\int_X \alpha\wedge S^\omega_1(\varphi)=0$. Hence $\int_X S^\omega_1(\varphi)\wedge \omega =0$, then we have $S^\omega_1(\varphi)\wedge \omega = 0$, which means
\begin{equation}\label{eq 3.3}
    \langle \omega+ dd^c \varphi \rangle \wedge \omega= (\omega+dd^c \varphi ) \wedge \omega.
\end{equation}

 We thus obtain the following result:
\begin{prop}\label{prop 3.1}
    Let $(X,\omega)$ be a compact Kähler surface and $\varphi \in DMA(X,\omega)\cap \mathcal{G}(X,\omega)$. Then,
\begin{align*}
{\rm MA}_\omega(\varphi)&=[\omega+dd^c \varphi]^2=\langle \omega+ dd^c \varphi \rangle\wedge (\omega + dd^c \varphi)\\
&= dd^c \varphi \wedge \langle \omega+ dd^c \varphi \rangle + \omega \wedge ( \omega + dd^c \varphi).
\end{align*}
In particular, when $\varphi \in \mathcal{D}(X,\omega)$, the above formula holds.
\end{prop}
Furthermore, for the class $\mathcal{D}(X,\omega)$, we have
\begin{thm}\label{thm 3.2}
    Let $\varphi_1,...,\varphi_n \in \mathcal{D}(X,\omega)$. Then, $\varphi_1 \in L^1\Big(\langle \omega_{\varphi_2} \wedge...\wedge \omega_{\varphi_n} \rangle\wedge\omega \Big)$.\par
    If $X$ is a compact Kähler surface and $\varphi,\psi \in \mathcal{D}(X,\omega)$, then we have
    $$
    {\rm MA}_\omega(\varphi,\psi)= \frac{1}{2}\left[ (\omega+dd^c\varphi)\wedge\langle \omega+dd^c \psi \rangle +(\omega+dd^c\psi)\wedge\langle \omega+dd^c \varphi \rangle \right].
    $$
\end{thm}
The proof needs the following Lemma:
\begin{lem}\label{lem 3.3}
    Let $\varphi$ be a quasi-psh function on $X$, and let $T$ be a positive closed $(n-1,n-1)$-current such that $\varphi \in L^1(\omega \wedge T)$. Then, for $\varphi_j \in C^\infty(X)$ such that $\varphi_j \searrow \varphi$, we have
    $$
    dd^c \varphi_j \wedge T \to dd^c \varphi \wedge T
    $$
    in the weak sense.
\end{lem}
%\begin{proof}   Without loss of generality, we can assume that $\varphi_1 \leq 0$. Set $f\in C^\infty(X)$ and $C\in \mathbb{R}^+$ such that $-C\omega \leq dd^c f\leq C\omega$. We have  \begin{align*}   \int_X f dd^c \varphi \wedge T &= \int_X \varphi dd^c f \wedge T\\  &= \int_X \varphi (C\omega + dd^c f)\wedge T - \int_X \varphi ~ C\omega \wedge T.   \end{align*}   By the monotone convergence theorem, we have   \begin{align*}   \int_X \varphi (C\omega + dd^c f)\wedge T &- \int_X \varphi ~ C\omega \wedge T\\  &=\lim_j \left( \int_X \varphi_j (C\omega + dd^c f)\wedge T - \int_X \varphi_j ~ C\omega \wedge T \right)\\  &=\lim_j \int_X \varphi_j dd^c f \wedge T =\lim_j \int_X f dd^c \varphi_j \wedge T,   \end{align*}   where in the last line we used Stokes theorem. \end{proof}

\begin{proof}[Proof of Theorem \ref{thm 3.2}]
     $1^\circ.$ 
     We first prove that $\varphi_1 \in L^1\Big(\langle \omega_{\varphi_2} \wedge...\wedge \omega_{\varphi_n} \rangle\wedge\omega \Big)$. Without loss of generality, we can assume that $\varphi_i \leq 0,~i=1,...,n$.\par
     It follows from Corollary \ref{cor 2.14} and Theorem \ref{thm 2.17} that
     $$
     \frac{1}{n}\left(\sum_i \varphi_i\right) \in L^1 \left(\langle (\omega+dd^c\frac{1}{n}(\sum_i \varphi_i))^{n-1}\rangle\wedge\omega \right).
     $$
     By \cite[Proposition 8.16]{GZ17}, there exists $\varphi^j_1\in {\rm PSH}(X,\omega) \cap C^\infty(X)$ such that $\varphi^j_1 \searrow \varphi_1$ as $j\to+\infty$.
     Then, by the multilinearity of the non-pluripolar product and multinomial theorem, we obtain
     $$
      -\int_X \varphi_1^j \omega \wedge \langle (\omega+dd^c \varphi_2) \wedge...\wedge (\omega+ dd^c \varphi_n) \rangle \leq \int_X \varphi_1 \langle (n\omega+dd^c(\sum_i \varphi_i))^{n-1}\rangle\wedge\omega.
     $$
     Applying the monotone convergence theorem, we obtain
     $$
     0\leq -\int_X \varphi_1 \omega \wedge \langle (\omega+dd^c \varphi_2) \wedge...\wedge (\omega+ dd^c \varphi_n) \rangle < +\infty.
     $$
     $2^\circ.$ Next, we prove that 
     $$
     {\rm MA}_\omega(\varphi,\psi)= \frac{1}{2}\left[ (\omega+dd^c\varphi)\wedge\langle \omega+dd^c \psi \rangle +(\omega+dd^c\psi)\wedge\langle \omega+dd^c \varphi \rangle \right].
     $$
     By Corollary \ref{cor 2.14} with $t_1=t_2=\frac{1}{2}$ and Proposition \ref{prop 3.1}, we obtain
    \begin{align*}
    {\rm MA}_\omega(\varphi,\psi) &= 2 {\rm MA}_\omega\left(\frac{1}{2}(\varphi+\psi ) \right)-\frac{1}{2}\left( {\rm MA}_\omega(\varphi)+{\rm MA}_\omega(\psi) \right)\\
    =&\frac{1}{2}\left( \omega\wedge \left(\omega +dd^c\varphi\right) + \omega\wedge \left(\omega +dd^c\psi\right) \right) \\
    &-\frac{1}{2} \left( dd^c \varphi \wedge \langle \omega+  dd^c\varphi  \rangle +  dd^c \psi \wedge \langle \omega+  dd^c\psi  \rangle \right)\\
    &+\left(\left( dd^c(\varphi+\psi)\right)\wedge \langle \omega+ 2^{-1} dd^c(\varphi+\psi)  \rangle \right).
    \end{align*}
    Therefore, we can complete the proof of the theorem if we can show that
    \begin{align*}
     &dd^c(\varphi+\psi)\wedge \langle \omega + 2^{-1} dd^c(\varphi+\psi)  \rangle  \\
    = \frac{1}{2} \Big( dd^c \varphi \wedge \langle \omega+  dd^c\varphi  \rangle& +  dd^c \psi \wedge \langle \omega +  dd^c\psi  \rangle +   dd^c \varphi \wedge \langle \omega+  dd^c\psi  \rangle +  dd^c \psi \wedge \langle \omega+  dd^c\varphi  \rangle \Big).
    \end{align*}
    Indeed, analogous to $1^\circ$, we can find smooth functions $\psi_j$ and $\varphi_j$ that decrease to $\psi$ and $\varphi$, respectively. Applying the multilinearity of the non-pluripolar product and Lemma \ref{lem 3.3}, makes it straightforward to obtain the desired conclusion.
\end{proof}

\section{Solving the complex Monge--Ampère equation}\label{Section 4}
\subsection{Cegrell--Lebesgue decomposition}
First, let us recall the Cegrell--Lebesgue decomposition theorem in a bounded hyperconvex domain $\Omega \subset \mathbb{C}^n$.
\begin{thm} \label{thm 4.1} $($\cite[Theorem 5.11]{Ce04}$)$
    Let $\mu$ be a positive measure on $\Omega$. Then, $\mu$ can be decomposed into a regular positive (non-pluripolar) measure $\mu_r$, and the singular positive (pluripolar) measure $\mu_s$ has support on some pluripolar set such that
    $$
    \mu= \mu_r + \mu_s.
    $$
\end{thm}
We claim that the Cegrell--Lebesgue decomposition of $\mu$ is unique. Indeed, suppose that $\mu_r'$ and $\mu_s'$ are other Cegrell--Lebesgue decompositions, and that 
$$
\mu_s= \mathbf{1}_{\{u=-\infty \}} \mu_s = \mathbf{1}_{\{u=-\infty \}} \mu,~~\mu_s' = \mathbf{1}_{ \{ u' = -\infty \} } \mu_s' = \mathbf{1}_{ \{ u' = -\infty \} }\mu,
$$
where $u,u' \in {\rm PSH}(\Omega) $. Set $A:= \{ u = -\infty\}\cup \{ u'=-\infty \} = \{ u+ u' = -\infty\}$. For $\forall f \in C_0 (\Omega)$, since $\int_A \mu_r = \int_A \mu_r'=0$, we have
$$
\int f d\mu_s = \int \mathbf{1}_A f d\mu_s = \int_A f d\mu = \int \mathbf{1}_A f d\mu_s' = \int f d\mu_s'.
$$
  This means that $\mu_s = \mu_s'$; hence, 
  $\mu_r = \mu - \mu_s = \mu- \mu_s ' = \mu_r '$.\par

  Let $(X,\omega_X)$ be a compact complex manifold of complex dimension $n$ equipped with a Hermitian form, and let $\mu$ be a positive measure on $X$. 
  Let $\{ U_i ,z_i \}_{i=1}^N$ be the holomorphic atlas of $X$ such that $z_i(U_i)$ is the unit ball of $\mathbb{C}^n$ and $\{ \chi_i \}_{i=1}^N$ is the partition of unity of $\{ U_i \}_i$. 
  According to Theorem \ref{thm 4.1}, we have $\chi_i \mu = \mu_r^i + \mu_s^i$ for each $i$, where $\mu_r^i$ is a positive non-pluripolar measure on $U_i$ and  $\mu_s^i$ is supported on the set $\{ u_i = -\infty \},~ u_i \in {\rm PSH}(U_i)$. 
  Set $\mu_r := \sum_i \mu_r^i$ and $\mu_s := \sum_i \mu_s^i$, we have the Cegrell--Lebesgue decomposition $\mu = \mu_r + \mu_s $ in the global setting.\par
  By the same method as above, the decomposition is unique in the global setting. Furthermore, thanks to \cite[Theorem 1.1]{Vu19}, there exists $u \in {\rm PSH}(X,\omega_X)$ such that $\mu_s$ is supported on the set $\{u = -\infty\}$.\par
  From the above discussion, we have the following result:
\begin{cor}\label{cor 4.4}
    Let $(X,\omega_X)$ be a compact complex manifold of complex dimension $n$ equipped
    with a Hermitian form. Let $\mu$ be a positive measure on $X$. Then, $\mu$ can be uniquely decomposed into a regular positive (non-pluripolar) measure $\mu_r$ and a singular positive (pluripolar) measure $\mu_s$ support on $\{ u = -\infty\},~u\in {\rm PSH } (X,\omega_X)$ such that
    $$
    \mu= \mu_r + \mu_s.
    $$
\end{cor}

  Now, we return to the compact Kähler manifold $(X,\omega)$. 
  From Corollary \ref{cor 4.4}, for a normalized measure $\mu$, a natural idea is to decompose $\mu$ into regular and singular parts $\mu_r$ and $\mu_s$ and consider
\begin{equation}\label{eq 4.1}
    {\rm MA}_\omega(\varphi)_r = \mu_r,~{\rm MA}_\omega(\varphi)_s = \mu_s,~\varphi \in \mathcal{D}(X,\omega),
\end{equation}
respectively. Fortunately, by \cite[Corollary 1.8]{CGZ08} we have:

\begin{thm}\label{thm 4.5}
    Let $\varphi \in \mathcal{D}(X,\omega)$. Then
    $$
    {\rm MA}_\omega(\varphi)_r = \langle \omega_\varphi^n \rangle~~{\rm and}~~{\rm MA}_\omega(\varphi)_s = \mathbf{1}_{ \{ \varphi = -\infty \}} {\rm MA}_\omega(\varphi)_s.
    $$
\end{thm}
Then, we can provide:
\begin{cor}\label{cor 4.6}
    Let $\varphi , \psi \in \mathcal{D}(X,\omega)$ such  that $\varphi \preceq \psi$. Then,
    $$
    {\rm MA}_\omega(\varphi)_s \geq {\rm MA}_\omega(\psi)_s.
    $$
\end{cor}
\begin{proof}
Let $\{ U_i, z_i \}_{i=1}^N$ be the holomorphic atlas of $X$ such that $z_i(U_i) \subset\mathbb{C}^n$ is the unit ball and $\omega  = dd^c g_i$ on $U_i$, where $g_i \in C^\infty \cap L^\infty\cap {\rm PSH}(U_i)$.
    In each $U_i$, we have $\varphi + g_i \leq \psi +g_i + C$ for some $C>0$. Then, by using \cite[Lemma 4.1]{ACCP09}, we have 
    \begin{equation}\label{eq 4.6}
    \begin{split}
        \mathbf{1}_{U_i}  \mathbf{1}_{ \{ \varphi +g_i=-\infty \}} & {\rm MA}(\varphi+g_i) \geq \mathbf{1}_{U_i}  \mathbf{1}_{ \{ \psi +g_i=-\infty \}}  {\rm MA}(\psi+g_i).
    \end{split}
    \end{equation}
    Let $\{ \chi_i \}_{i=1}^N$ be the partition of unity of $\{ U_i \}$. Since ${\rm MA}_\omega (\varphi) = \sum_i \chi_i \mathbf{1}_{U_i}  {\rm MA}(g_i+\varphi)$, by (\ref{eq 4.6}), we are done.
\end{proof}

When $\varphi \in \mathcal{D}(X,\omega)$, due to the uniqueness of the Cegrell–Lebesgue decomposition, equation (\ref{eq 3.1}) and Theorem \ref{thm 4.5} yield the following result.
\begin{cor}
    Let $\varphi \in \mathcal{D}(X,\omega)$. Then,
    $$
    {\rm MA}_\omega(\varphi)_r = \langle  \omega_\varphi^n \rangle~~~~{\rm and} ~~~~
    {\rm MA}_\omega(\varphi)_s = \sum_{j=1}^n S^\omega_j(\varphi) \wedge \omega^{n-j}.
    $$
\end{cor}

\subsection{Solving the complex Monge--Ampère equation in $\mathcal{D}(X,\omega)$}\label{sec 4.2}
In this subsection, we assume $(X,\omega)$ is a Kähler manifold such that $\int_X \omega^n=1$.\par
We set a normalized measure $\mu = \mu_r +\mu_s$ on $X$ such that $\int_X \mu=1$. If there exists $\varphi \in \mathcal{D}(X,\omega)$ such that ${\rm MA}_\omega (\varphi)=\mu$, it follows from Corollary \ref{cor 4.4} and Theorem \ref{thm 4.5} that
$$
\int_X \langle \omega_\varphi^n \rangle = \int_X \mu_r ~~{\rm and}~~ 
\mathbf{1}_{ \{ \varphi=-\infty\} } {\rm MA}_\omega(\varphi) = \mu_s.
$$
Then we have
\begin{prop}\label{prop 4.8}
    Assume that $\varphi \in \mathcal{D}(X,\omega)$ and $\mu_s = {\rm MA}_\omega(\varphi)_s$. Then, for $\psi \in {\rm PSH}(X,\omega)$ so that $\varphi \preceq \psi \preceq P_\omega[\varphi]$, we have ${\rm MA}(\psi)_s = \mu_s$. Furthermore, if $\psi \succeq \varphi$ so that ${\rm MA}(\psi)_s = \mu_s$, then we have $\psi \preceq P_\omega[\varphi]$.
\end{prop}
\begin{proof}
 Assume that $\psi \in {\rm PSH}(X,\omega)$ such that $\varphi \preceq \psi \preceq P_\omega[\varphi]$. By Proposition \ref{thm 2.2} and Lemma \ref{lem 2.11}, we have $\psi \in \mathcal{E}(X,\omega,P_\omega[\varphi])\cap \mathcal{D}(X,\omega)$. Applying Corollary \ref{cor 4.6}, we obtain
\begin{align*}
{\rm MA}_\omega(\varphi)_s \geq  {\rm MA}_\omega (\psi)_s.
\end{align*}
Note that we have
\begin{align*}
    \int_X {\rm MA}_\omega(\psi)_s &= 1 - \int_X
    {\rm MA}_\omega(\psi)_r = 1 - \int_X \langle \omega_\psi^n \rangle\\
    &= 1 - \int_X \langle \omega_\varphi^n \rangle = 1- \int_X 
    {\rm MA}_\omega(\varphi)_r =  \int_X {\rm MA}_\omega(\varphi)_s,
\end{align*}
where we used the definition of $\mathcal{E}(X,\omega,P_\omega[\varphi])$ and Theorem \ref{thm 4.5}.  By comparing the total mass, we obtain 
$\mu_s={\rm MA}_\omega(\varphi)_s={\rm MA}_\omega(P_\omega[\varphi])_s$.\par
If $\psi \succeq \varphi$ so that ${\rm MA}(\psi)_s = {\rm MA}(\varphi)_s$, then we have $\int_X \langle \omega_\psi^n \rangle = \int_X \langle \omega_\varphi^n \rangle$. It follows from Proposition \ref{DDL23 2.4} (ii). that $\psi \preceq P_\omega[\varphi]$.
\end{proof}
 Set $\nu  = \nu_r + \mu_s$. If $\nu_r $ satisfies $\int_X \nu_r = 1 - \int_X \mu_s$ and $\nu = f \omega^n$ for some $f \in L^p(\omega),~p>1$, it follows from \cite[Theorem 5.20]{DDL23} that there exists $\psi \in {\rm PSH}(X,\omega)$ such that 
$$
\psi \simeq P_\omega[\varphi]~~{\rm and}~~{\rm MA}_\omega(\psi)_r = \langle  \omega_\psi^n \rangle = \nu_r.
$$
By Lemma \ref{lem 2.11} and Proposition \ref{prop 4.8}, we have $\psi \in \mathcal{D}(X,\omega)$ and ${\rm MA}_\omega(\psi)_s = \mu_s$. We thus obtain
\begin{align*}
{\rm MA}_\omega(\psi) &= {\rm MA}_\omega(\psi)_r +{\rm MA}_\omega(\psi)_s= \langle \omega_\psi^n \rangle + \mu_s = \nu_r +\mu_s.
\end{align*}

 In summary, we obtain the following Theorem.
\begin{thm}\label{thm 4.9}
    Let $(X,\omega)$ be a compact Kähler manifold of complex dimension $n$ such that $\int_X \omega^n =1$, and let $\mu_s$ be a pluripolar measure on $X$ supported on some pluripolar set. If there exists $\varphi \in \mathcal{D}(X,\omega)$ such that 
    ${\rm MA}_\omega(\varphi)_s = \mu_s$, then for all $0 \leq f \in L^p(\omega^n),~p>1$ such that $\int_X f \omega^n = 1 - \int_X \mu_s$, there exists a $\psi \in {\rm PSH}(X,\omega)$ that is the solution of the equation
    $$
    {\rm MA}_\omega(\psi)= f \omega^n + \mu_s,~\psi \in \mathcal{D}(X,\omega).
    $$
\end{thm}

\subsection{Application to compact Kähler surfaces} 
  Let $(X,\omega)$ be a compact Kähler surface such that $\int_X \omega^2=1$, and let $\mu$ be a measure on $X$ so that $\mu=\mu_s$. 
  If there exists $\varphi \in \mathcal{D}(X,\omega)$ such that ${\rm MA}_\omega(\varphi) = \mu$ (e.g. the quasi-psh Green functions, see \cite{CG09}), then, for $\forall s \in [0,1]$, it follows from Corollary \ref{cor 2.14}, (\ref{eq 3.3}), and Theorem \ref{thm 3.2} that
    \begin{align*}
       {\rm MA}_\omega(s \varphi) &= (1-s)^2 {\rm MA}_\omega(0) + 2(1-s)s
       {\rm MA}_\omega (0,\varphi)+ s^2 {\rm MA}_\omega(\varphi)\\
       &= (1-s)^2 \omega^2 + 2(1-s)s \omega\wedge \langle \omega+ dd^c \varphi \rangle + s^2 \mu_s.
    \end{align*}
   This means that $ {\rm MA}_\omega(s \varphi)_s = s^2 \mu$.
   Therefore, for any $t \in [0,1]$ and any $0\leq f \in L^p(\omega^n)$, $p>1$ such that $\int_X f\omega^2 = 1$, by Theorem \ref{thm 4.9}, there exists $\psi_{f,t} \in \mathcal{D}(X,\omega)$ that satisfies
   $$
     {\rm MA}_\omega(\psi_{f,t}) = (1-t) f\omega^2 + t\mu.
   $$

\section{Acknowledgments} I want to thank my Master's advisor Liyou Zhang for interesting discussions and valuable suggestions on an early draft, which made the paper more readable. I would also like to thank T. Darvas for his insightful comments that improved the presentation of the paper.

\end{document}